\documentclass[12pt,twoside,leqno]{article}
\usepackage{amsmath}
\usepackage{amssymb}
\usepackage{amsxtra}
\usepackage{amscd}
\usepackage{amsthm}
\usepackage[mathscr]{eucal}

\theoremstyle{plain}

\newtheorem{lem}[subsection]{Lemma}
\newtheorem{conj}[subsection]{Conjecture}

\theoremstyle{definition}

\newtheorem{para}[subsection]{}

\newenvironment{pf}{\proof[\proofname]}{\endproof}

\begin{document}

\title{Heights of mixed motives}

\author{Kazuya Kato}

\maketitle

\newcommand\Cal{\mathcal}
\newcommand\define{\newcommand}

\define\gp{\mathrm{gp}}%
\define\fs{\mathrm{fs}}%
\define\an{\mathrm{an}}%
\define\mult{\mathrm{mult}}%
\define\Ker{\mathrm{Ker}\,}%
\define\Coker{\mathrm{Coker}\,}%
\define\Hom{\mathrm{Hom}\,}%
\define\Ext{\mathrm{Ext}\,}%
\define\rank{\mathrm{rank}\,}%
\define\gr{\mathrm{gr}}%
\define\cHom{\Cal{Hom}}
\define\cExt{\Cal Ext\,}%

\define\cA{\Cal A}
\define\cC{\Cal C}
\define\cD{\Cal D}
\define\cO{\Cal O}
\define\cS{\Cal S}
\define\cM{\Cal M}
\define\cG{\Cal G}
\define\cH{\Cal H}
\define\cE{\Cal E}
\define\cF{\Cal F}
\define\cN{\Cal N}
\define\fF{\frak F}
\define\Dc{\check{D}}
\define\Ec{\check{E}}

\newcommand{\A}{{\mathbb{A}}}
\newcommand{\N}{{\mathbb{N}}}
\newcommand{\Q}{{\mathbb{Q}}}
\newcommand{\Z}{{\mathbb{Z}}}
\newcommand{\R}{{\mathbb{R}}}
\newcommand{\C}{{\mathbb{C}}}
\newcommand{\bN}{{\mathbb{N}}}
\newcommand{\bQ}{{\mathbb{Q}}}
\newcommand{\bF}{{\mathbb{F}}}
\newcommand{\bZ}{{\mathbb{Z}}}
\newcommand{\bP}{{\mathbb{P}}}
\newcommand{\bR}{{\mathbb{R}}}
\newcommand{\bC}{{\mathbb{C}}}
\newcommand{\bbQ}{{\bar \mathbb{Q}}}
\newcommand{\ol}[1]{\overline{#1}}
\newcommand{\too}{\longrightarrow}
\newcommand{\respect}{\rightsquigarrow}
\newcommand{\compatible}{\leftrightsquigarrow}
\newcommand{\upc}[1]{\overset {\lower 0.3ex \hbox{${\;}_{\circ}$}}{#1}}
\newcommand{\Gmlog}{\bG_{m, \log}}
\newcommand{\Gm}{\bG_m}
\newcommand{\ep}{\varepsilon}
\newcommand{\Spec}{\operatorname{Spec}}
\newcommand{\val}{{\mathrm{val}}} 
\newcommand{\n}{\operatorname{naive}}
\newcommand{\bs}{\operatorname{\backslash}}
\newcommand{\Gal}{\operatorname{{Gal}}}
\newcommand{\gal}{{\rm {Gal}}({\bar \Q}/{\Q})}
\newcommand{\galp}{{\rm {Gal}}({\bar \Q}_p/{\Q}_p)}
\newcommand{\gall}{{\rm{Gal}}({\bar \Q}_\ell/\Q_\ell)}
\newcommand{\wep}{W({\bar \Q}_p/\Q_p)}
\newcommand{\wel}{W({\bar \Q}_\ell/\Q_\ell)}
\newcommand{\Ad}{{\rm{Ad}}}
\newcommand{\BS}{{\rm {BS}}}
\newcommand{\even}{\operatorname{even}}
\newcommand{\End}{{\rm {End}}}
\newcommand{\odd}{\operatorname{odd}}
\newcommand{\GL}{\operatorname{GL}}
\newcommand{\np}{\text{non-$p$}}
\newcommand{\g}{{\gamma}}
\newcommand{\G}{{\Gamma}}
\newcommand{\Lam}{{\Lambda}}
\newcommand{\La}{{\Lambda}}
\newcommand{\lam}{{\lambda}}
\newcommand{\la}{{\lambda}}
\newcommand{\uL}{{{\hat {L}}^{\rm {ur}}}}
\newcommand{\uQp}{{{\hat \Q}_p}^{\text{ur}}}
\newcommand{\sel}{\operatorname{Sel}}
\newcommand{\dt}{{\rm{Det}}}
\newcommand{\Sig}{\Sigma}
\newcommand{\fil}{{\rm{fil}}}
\newcommand{\SL}{{\rm{SL}}}
\newcommand{\spl}{{\rm{spl}}}
\newcommand{\st}{{\rm{st}}}
\newcommand{\Isom}{{\rm {Isom}}}
\newcommand{\Mor}{{\rm {Mor}}}
\newcommand{\bg}{\bar{g}}
\newcommand{\id}{{\rm {id}}}
\newcommand{\cone}{{\rm {cone}}}
\newcommand{\al}{a}
\newcommand{\ChL}{{\cal{C}}(\La)}
\newcommand{\Image}{{\rm {Image}}}
\newcommand{\toric}{{\operatorname{toric}}}
\newcommand{\torus}{{\operatorname{torus}}}
\newcommand{\Aut}{{\rm {Aut}}}
\newcommand{\Qp}{{\mathbb{Q}}_p}
\newcommand{\barQp}{{\mathbb{Q}}_p}
\newcommand{\Qpur}{{\mathbb{Q}}_p^{\rm {ur}}}
\newcommand{\Zp}{{\mathbb{Z}}_p}
\newcommand{\Zl}{{\mathbb{Z}}_l}
\newcommand{\Ql}{{\mathbb{Q}}_l}
\newcommand{\Qlur}{{\mathbb{Q}}_l^{\rm {ur}}}
\newcommand{\F}{{\mathbb{F}}}
\newcommand{\eps}{{\epsilon}}
\newcommand{\epsLa}{{\epsilon}_{\La}}
\newcommand{\epsLaVxi}{{\epsilon}_{\La}(V, \xi)}
\newcommand{\epsOLaVxi}{{\epsilon}_{0,\La}(V, \xi)}
\newcommand{\Qplin}{{\mathbb{Q}}_p(\mu_{l^{\infty}})}
\newcommand{\otimesQplin}{\otimes_{\Qp}{\mathbb{Q}}_p(\mu_{l^{\infty}})}
\newcommand{\galFl}{{\rm{Gal}}({\bar {\Bbb F}}_\ell/{\Bbb F}_\ell)}
\newcommand{\gallur}{{\rm{Gal}}({\bar \Q}_\ell/\Q_\ell^{\rm {ur}})}
\newcommand{\galFF}{{\rm {Gal}}(F_{\infty}/F)}
\newcommand{\galFv}{{\rm {Gal}}(\bar{F}_v/F_v)}
\newcommand{\galF}{{\rm {Gal}}(\bar{F}/F)}
\newcommand{\epsVxi}{{\epsilon}(V, \xi)}
\newcommand{\epsOVxi}{{\epsilon}_0(V, \xi)}
\newcommand{\plim}{\lim_
{\scriptstyle 
\longleftarrow \atop \scriptstyle n}}
\newcommand{\sig}{{\sigma}}
\newcommand{\ga}{{\gamma}}
\newcommand{\del}{{\delta}}
\newcommand{\Vss}{V^{\rm {ss}}}
\newcommand{\Bst}{B_{\rm {st}}}
\newcommand{\Dpst}{D_{\rm {pst}}}
\newcommand{\Dcrys}{D_{\rm {crys}}}
\newcommand{\DdR}{D_{\rm {dR}}}
\newcommand{\Fin}{F_{\infty}}
\newcommand{\Kla}{K_{\lambda}}
\newcommand{\Ola}{O_{\lambda}}
\newcommand{\Mla}{M_{\lambda}}
\newcommand{\Det}{{\rm{Det}}}
\newcommand{\Sym}{{\rm{Sym}}}
\newcommand{\LaSa}{{\La_{S^*}}}
\newcommand{\cX}{{\cal {X}}}
\newcommand{\MHG}{{\frak {M}}_H(G)}
\newcommand{\tauMla}{\tau(M_{\lambda})}
\newcommand{\Fvur}{{F_v^{\rm {ur}}}}
\newcommand{\Lie}{{\rm {Lie}}}
\newcommand{\cL}{{\cal {L}}}
\newcommand{\cW}{{\cal {W}}}
\newcommand{\fq}{{\frak {q}}}
\newcommand{\cont}{{\rm {cont}}}
\newcommand{\SC}{{SC}}
\newcommand{\Om}{{\Omega}}
\newcommand{\dR}{{\rm {dR}}}
\newcommand{\crys}{{\rm {crys}}}
\newcommand{\hatSig}{{\hat{\Sigma}}}
\newcommand{\rdet}{{{\rm {det}}}}
\newcommand{\ord}{{{\rm {ord}}}}
\newcommand{\BdR}{{B_{\rm {dR}}}}
\newcommand{\BdRO}{{B^0_{\rm {dR}}}}
\newcommand{\Bcrys}{{B_{\rm {crys}}}}
\newcommand{\Qw}{{\mathbb{Q}}_w}
\newcommand{\barkappa}{{\bar{\kappa}}}
\newcommand{\cP}{{\Cal {P}}}
\newcommand{\cZ}{{\Cal {Z}}}
\newcommand{\oppLa}{{\Lambda^{\circ}}}

\begin{abstract}
 We define the height of a mixed motive over a number field extending our previous work for pure motives.  
\end{abstract}
\renewcommand{\thefootnote}{\fnsymbol{footnote}}
\footnote[0]{Primary 14G40; Secondary 14G25,11G50. 
The author is partially supported by an NSF grant.}

In \cite{KK}, we defined the height of a pure motive over a number field. In this paper, we extend the definition to mixed motives.

For a mixed motive $M$ over a number field,  the height $h(M)$ is defined as the sum of the heights $h_{w,d}(M)$ for $w\in \Z$ and $d\geq 0$, and to define $h_{w,d}$ for $d\geq 1$, we have to fix a polarization of the pure motive $P_{w,k}:=(\gr^W_wM)^*\otimes \gr^W_{w-d}M$ where $W$ is the weight filtration of $M$. 

In the case $d=0$, $h_{w,0}(M)$ is the  height of the pure motive $h(\gr^W_wM)$ defined in \cite{KK}.

The case $d=1$ is 
the height of of  the extension $0\to \gr^W_{w-1}M\to W_wM/W_{w-2}M\to \gr^W_wM\to 0$ defined by Beilinson \cite{Be} and Bloch  \cite{Bl}.  We review the idea of the definition in \ref{d=1}. 

For  $d\geq 2$,  $h_{w,d}(M) $ is  determined by $W_wM/W_{w-d-1}M$. It is defined in section \ref{sec.ari} as the sum of local heights.
We first give the geometric analogue of $h_{w,d}$ ($d\geq 2$)  in section \ref{sec.geo}, for the idea in the number field case can be seen well by comparison.

Ideas in this paper came from a joint work of Spencer Bloch  and the author. Also, joint works with Chikara Nakayama and Sampei Usui inspired this work. 

Details will be published elsewhere.
\section{Review on the splitting of Deligne}\label{prim}

This  section is a preliminary for sections \ref{sec.geo} and \ref{sec.ari}. 

\begin{para}

Let $V$ be an abelian group, let $W=(W_w)_{w\in \bZ}$ be a finite increasing filtration on $V$, and let $N:V\to V$ be a nilpotent homomorphism such that $NW_w\subset W_w$ for alll $w\in \Z$.  

Then a finite increasing filtration $W'=(W'_w)_{w\in \bZ}$ on $V$ is called the 
{\it relative monodromy filtration of $N$ with respect to $W$} if it satisfies the following conditions (i) and (ii). 
(i) $NW'_w\subset W'_{w-2}$ for any $w\in \bZ$.  (ii) For any $w\in \bZ$ and $m\geq 0$, we have an isomorphism $N^m: \gr^{W'}_{w+m}\gr^W_w\overset{\cong}\longrightarrow  \gr^{W'}_{w-m}\gr^W_w.$

The relative monodromy filtration of $N$ with respect to $W$ need not exist. If it exists, it is unique (\cite{De} 1.6.13).
In the rest of section \ref{prim}, we assume that the relative monodromy filtration $W'$ of $N$ with respect to $W$ exists. 
\end{para}

\begin{para}
Let $m\geq 0$. Then $\gr^{W'}_{w+m}\gr^W_w=A \oplus B$, where $A$ is the kernel of $\gr^{W'}_{w+m}\gr^W_w\overset{N^{m+1}}\longrightarrow \gr^{W'}_{w-m-2}\gr^W_w$ and $B$ is the image of $N:\gr^{W'}_{w+m+2}\gr^W_w\to \gr^{W'}_{w+m}\gr^W_w$. The component $A$ is called the {\it primitive component} of $\gr^{W'}_{w+m}\gr^W_w$ and will be denoted by $(\gr^{W'}_{w+m}\gr^W_w)_{\text{prim}}$.

\end{para}

\begin{para}\label{primN1} Denote by $W_{\bullet}Hom(V, V)$ (resp. $W'_{\bullet}Hom(V, V)$)
the filtration on $Hom(V, V)$ induced by $W$ (resp. $W'$). Then $W'_{\bullet}Hom(V, V)$ is the relative weight filtration of the nilpotent homomorphism $Ad(N): Hom(V, V)\to Hom(V, V)$ with respect to $W_{\bullet}Hom(V, V)$. 

\end{para}

\begin{para}\label{Dspl} Assume  $V=\oplus_w \; U'_w$ is a splitting of  $W'$ (that is, $W'= \sum_{n\leq w} U'_n$ for any $w$) which satisfies  $N(U'_w)\subset U'_{w-2}$ for all $w$ and which is compatible with $W$ (that is, $W_w= \sum_n (W_w\cap U'_n)$ for any $n$).
Then by Deligne (see \cite{Sc}),   there is a unique splitting $V= \oplus_w U_w$ of $W$ such that $V=\oplus_{m,n} U_m\cap U'_n$ and such that 
  if  $N_w\in \gr^W_wHom(V, V)$ denotes the component of weight $w$ of $N$ with respect to this splitting of $W$ and if $\bar N_w$ denotes the class
of $N_w$ in $\gr^{W'}_{-2}\gr^W_wHom(V, V)$, then $\bar N_{-1}=0$ and for any $w\leq -2$, $\bar N_w$  belongs to the primitive component of $\gr^{W'}_{-2}\gr^W_wHom(V,V)$.

\end{para}

\begin{lem}\label{primN} Assume that a splitting of $W'$ as in \ref{Dspl} exists. Then the element
$\bar N_w$ of $(\gr^{W'}_{-2}\gr^W_wHom(V, V))_{\text{prim}}$ is independent of the choice of such splitting of $W'$. 

\end{lem}

\begin{pf}  Assume we are given two such splittings $((U'_w)^{(i)})$ ($i=1,2$) of $W'$, and let $g$ be the unique automorphism of $(V,W')$ which induces the identity map on $\gr^{W'}$ such that $g(U'_w)^{(1)}=(U'_w)^{(2)}$. Then $gW=W$ and $gN=Ng$. For Deligne's splitting $(U^{(i)}_w)$ of $W$ associated to $((U'_w)^{(i)})$, we have $U^{(2)}_w=gU^{(1)}_w$. If $N_w^{(i)}$ and $\bar N_w^{(i)}$ denote the $N_w$ and $\bar N_w$ for $(U_w^{(i)})$, respectively, we have $N^{(2)}_w=gN^{(1)}_wg^{-1}$ and hence $\bar N^{(2)}_w= g\bar N_w^{(1)}g^{-1}=\bar N_w^{(1)}$.
\end{pf}

\section{Geometric heights}\label{sec.geo}

\begin{para}
We consider the case of variation of MHS (mixed Hodge structure) first.

Let $C$ be a proper smooth curve over $\C$ and let $H=(H_{\Z}, W, F)$ be a variation of MHS  on $C-S$ for some finite subset $S$ of $C$. Assume that 
the graded quotients $\gr^W_wH$ are polarized and assume that $H$ is admissible at 
any  point of $S$. 

Here the admissibility means that IHMH (infinitesimal mixed Hodge module) in the sense of Kashiwara \cite{KM} appears at each point of $C$.  In other words, the admissibility means that the following conditions (i) and (ii) are satisfied at each $x\in C$. (i) Let $N_x: H_{\Q,x}\to H_{\Q,,x}$ be the logarithm of the local monodormy at $x$. Then 
the relative monodromy filtration of $N_x$ with respect to $W_x$ exists.  (ii) The limit Hodge filtration appears at $x$. It is known that any variation of MHS with geometric origin satisfies these conditions (i) and (ii).

 Let $w, d\in \Z$ and assume $d\geq 2$.

 We will define the height $h_{w,d}(H)$ of $H$ as the sum of local heights $h_{w,d,x}(H)$ for all points $x$ of $C$. The local height $h_{w,d,x}(H)$ will be defined as the ''size'' of  $N_x$.

\end{para}
\begin{para}\label{dspl}  
 Let $W'$ be the relative monodromy filtration of $N_x$ with respect to $W_x$. By Kashiwara \cite{KM}, there is a splitting of $W'$ as in \ref{Dspl}. Hence by \ref{primN}, 
 $\bar N_{x,-d}\in (\gr^{W'}_{-2}\gr^W_{-d}Hom(H_x, H_x))_{\text{prim}}$ is defined for $d\geq 2$.
\end{para}

\begin{para}\label{pos} Now we define $h_{w,d,x}(H)$. Let $P=(\gr^W_wH)^*\otimes \gr^W_{w-d}H$. We have the component $\bar N_{x,w,d} \in \gr^{W'}_{-2}P$ of $\bar N_{x,-d}$. Let $\langle\;,\;\rangle: P \times P\to \Q$ be the pairing defined by the given polarizations of $\gr^W_wH$ and $\gr^W_{w-d}H$. Then by \cite{CKS}, \cite{KM}, the induced pairing
$$\langle\;,\;\rangle_{d-2}: \gr^{W'}_{-2}(P )\times \gr^{W'}_{-2}(P) \to \Q\;;\; (u,v) \mapsto \langle Ad(N_x)^{d-2}(u), v\rangle $$ satisfies 
$\langle \bar N_{x,w,d}, \bar N_{x,w,d}\rangle_{d-2}\geq 0$.

 We define  the local geometric height $h_{w,d,x}(H)$ by 
$$ h_{w,d,x}(H)= (\langle \bar N_{x,w,d}, \bar N_{x,w,d}\rangle_{d-2})^{1/d}\in \R_{\geq 0}.$$

(It is important to take the $d$-th root here, to have formulas like \ref{GA1}, \ref{GA2}.) 

\end{para}

\begin{para}
We define the global geometric height $h_{w,d}(M)$  by
$$h_{w,d}(H) = \sum_{x\in C}\; h_{w,d,x}(H).$$
\end{para}

\begin{para}\label{geoKk}  Let $K$ be a function field in one variable over a field $k$ and let $M$ be a mixed motive over $K$.

If $k$ is of characteristic $0$,  the definition of the geometric height $h_{w,d}(M)$ ($d\geq 2$) as well as that of  the local version  $h_{w,d,x}(M)$  goes in the same way as above, by using  the $\ell$-adic realization of $M$ in place of the local system $H_{\Q}$ of MHS. The fact it works independently of $\ell$ is reduced to the above story of MHS, by taking a subfield $k'$ of $k$ which is finitely generated over $\Q$ such that  $K$ and $M$ are defined over $k'$, and by embedding $k'$ into $\C$ and by using the variation of MHS associated to $M$.
\end{para}

\begin{para} \label{geop}

In the case $k$ is of characteristic $p>0$,  for the $\ell$-adic local system ($\ell\neq p$) associated to $M$, we have the following (WM) by Deligne \cite{De}. 

\medskip

(WM) The relative weight filtration $W'$ of $N_x$ with respect to $W_x$ exists, and it coincides with the weight filtration in the sense of mixed sheaf.

\medskip

 Since $W'$ has a splitting as in \ref{Dspl} by eigen values of a frobenius,  $\bar N_{x,d}$ ($d\geq 2$) are defined. 
The independence of $\ell$ and the positivity in \ref{pos}   may become questions.

\end{para}

\section{Heights of mixed motives over number fields}\label{sec.ari}

\begin{para}
Now we consider a mixed motive $M$ over a number field $K$ with polarized graded quotients for the weight filtration. The height $h_{w,d}(M)$ ($d\geq 2$) will be defined as the sum of local heights $h_{w,d,v}(M)$ for all places $v$ of $K$. 
\end{para}

\begin{para}
For a finite place $v$ of $K$, we assume that the definition in section \ref{sec.geo}
works by considering the associated $\ell$-adic Galois representation $M_{\ell}$ for $\ell$ not equal to the characteristic of the residue field of $v$, and (in place of $N_x$ in section \ref{sec.geo}) the logarithm of the action of a generator of the tame inertia group at $v$. 
We define
$h_{w,d,v}(M)$ to be $\log(N(v))$ times the local height defined by the method of section \ref{sec.geo}.

(A difficulty comes from the fact the analogue of (WM) in \ref{geop}  (the weight-monodromy conjecture)  is not yet proved. The independence of $\ell$ and the positivity in  \ref{pos} are  also problems. We can say only that we have a good definition  in many cases. )
\end{para}

\begin{para} Now we consider an Archimedean place $v$. Then $h_{w,d,v}(M)$ 
 is defined by using the Hodge structure associated to $M$ at $v$ as follows.
\end{para}

\begin{para} Let $H$ be a mixed Hodge structure. We review a homomorphism $\delta:\gr^W_w H_\R \to \gr^W_{w-d}H_\R$. 

As is explained in \cite{CKS} 2.20,  there is a unique pair $(\delta, \tilde F)$ where $\tilde F$ is a decreasing filtration on $H_{\C}$ such that $(H_{\R}, \tilde F)$ is an $\R$-split mixed Hodge structure, $\delta$ is a nilpotent $\R$-linear map $H_\R\to H_\R$ such that the original Hodge filtration $F$ of $H$ is expressed as $F=\exp(i\delta)\tilde F$ and such that the Hodge $(p,q)$-component $\delta_{p,q}$ of $\delta$ with respect to $\tilde F$ is zero unless $p<0$ and $q<0$. 

Furthermore,  a nilpotent linear map $\zeta:H_\R\to H_\R$ is defined by a universal Lie polynomial in $(\delta_{p,q})_{p,q}$ as in [CKS] 6.60.  
The $\R$-split mixed Hodge structure $(H_\R, \hat F)$ with $\hat F=\exp(\zeta)\tilde F$ and the $\R$-splitting of $W$ associated to $\hat F$ are  important (see \cite{CKS}, \cite{KNU}).

Let $\delta_{w,d}: \gr^W_wH_\R\to \gr^W_{w-d}H_\R$ be the linear map induced by the component of $\delta$ of weight $-d$ with respect to the last $\R$-splitting of $W$.
\end{para}

\begin{para} Let $H$ be a MHS. 
Let $P=(\gr^W_wH)^*\otimes \gr^W_{w-d}H$ and assume that $P$ is endowed with a polarization $p$. Then we define the height $h_{w,d}(H, p)$ 
as 
$$h_{w,d}(H, p):= ((\delta_{w,d}, \delta_{w,d})_p)^{1/d}$$
where $(\;,\;)_p$ is the Hodge metric (Hermitian form) associated to $p$.

\end{para}

\begin{para}

Let $M$ be a mixed motive over a number field $K$ with polarized graded quotient for the weight filtration.

For a real (resp. complex) place $v$ of $K$, we define the local height $h_{w,d,v}(M)$ as $ h_{w,d}(H,p)$ (resp. $2\cdot h_{w,d}(H,p)$) where $H$ is the mixed Hodge structure associated to $M$ at $v$ and $p$ is the induced polarization on $(\gr^W_wH)^*\otimes \gr^W_{w-d}H$. 

We define the global height $h_{w,d}(M)$ by
$$h_{w,d}(M):=\sum_v h_{w,d,v}(M)$$
where $v$ ranges over all places of $K$.
\end{para}

\begin{para} Examples.

1. Let $a\in K^\times$. Then $a$ corresponds to a mixed $\Z$-motive $M$ over $K$ with an exact sequence 
$0\to \Z(1)\to M\to \Z\to 0$. We have $h_{0,2,v}(M) =|\log(|a|_v)|$ for any $v$.

2. Let $E$ be an elliptic curve over $K$, consider the pure $\Z$-motive  $H^1(E)(2)$ over $K$ of weight $-3$, and consider a mixed $\Z$-motive $M$ over $K$ with an exact sequence $0\to H^1(E)(2)\to M\to \Z\to 0$. Then $M$ gives an element $a$ of $K_2(E)\otimes \Q$, and  $h_{0,3,v}(E)$ is expressed by the $K_2$-regulator  (resp. the tame symbol) of $a$ at $v$ if $v$ is Archimedean (resp. finite).

\end{para}
\section{On the cases $d=0, 1$}\label{d=01}

We describe the geometric analogues of $h_{w,0}$ (\ref{geo0},  \ref{geop}), and describe shortly the idea of the definition of $h_{w,1}$ by Beilinson and Bloch (\ref{d=1}). 

We define $h_{w,0}$ of a mixed object as the height  of the pure object $\gr^W_w$. 

\begin{para}\label{geo0} Let $(C, H)$ be as in section \ref{sec.geo} and assume $H$ is pure.  We define the height $h(H)$ of $H$ by
$$h(H) = \sum_{r\in \Z} \; r\cdot deg (\gr^r H_{\cO})$$
where $H_{\cO}$ is the vector bundle on $C$ corresponding to $H$ and $\gr^r=F^r/F^{r+1}$ for the Hodge filtration $F$ on $H_{\cO}$. Here the degree of a vector bundle means the degree of its highest exterior power.

This invariant is related to the works \cite{Gri} section 7, \cite{Pe}, \cite{Pe2}, \cite{VZ}, etc. The author is thankful to T. Koshikawa for pointing out this relation. Comparing with these works on the geometric analogue, he suggests that any constant $c<4/(\sum_{p,q} (p-q)^2h^{p,q})$ may work as $c$ in Conjecture 4.8 of \cite{KK}. 

The case of pure motives over $k(C)$ for $k$ of characteristic $0$ is similar. 

\end{para}

\begin{para}\label{geop}
For a proper smooth curve $C$ over a field $k$ of characteristic $p>0$, and for a $p$-adically integral  $F$-crystal $D$ (with $F: \varphi^*D\to D$) on $C$ with logarithmic poles  associated to a pure motive, its height should be defined to be
$$p^{-1}deg(D/F(\varphi^*D))$$  This will be discussed in our forthcoming paper.

\end{para}

\begin{para}\label{d=1} 
We review the idea of the definition of $h_{w,1}(M)$ by Beilinson and Bloch, shortly. Here to give a short explanation, we assume that usual philosophies on mixed motives are true. Let $P=(\gr^W_wM)^*\otimes \gr^W_{w-1}M$. The extension $0\to \gr^W_{w-1}M\to W_wM/W_{w-2}\to \gr^W_wM\to 0$ corresponds to an extension $0\to P\to Q \to \Z\to 0$. By taking the dual, we have an extension $0\to \Z(1) \to Q^*(1) \to P^*(1)\to 0$. Assume that $P$ is polarized. By pulling back by the polarization $P\to P^*(1)$, the last extension gives an extension $0\to \Z(1) \to Q' \to P\to 0$. By the exact sequence $Ext^1(\Z, Q')\to Ext^1(\Z, P)\to Ext^2(\Z, \Z(1))$ and by $Ext^2(\Z, \Z(1))=0$, we see that there is an element $a$ of $Ext^1(\Z, Q')$ whose image in $Ext^1(\Z, P)$ coincides with the class of $Q$. Then $a$ gives 
a mixed motive $R$ such that $W_0R=R$, $\gr^W_0R=\Z$, $\gr^W_{-1}R=P$, $\gr^W_{-2}R=\Z(1)$, and $W_{-3}R=0$. Define
$$h_{w,1}(M):= h_{0,2}(R).$$ This does not depend on the choice of $a\in Ext^1(\Z, Q')$ as above. 

The geometric analogue of $h_{w,1}$ is given similarly. 

\end{para}

\section{Some topics}

\begin{para}
For a mixed $\Z$-motive $M$ over a number field $K$ whose graded quotients for the weight filtration are polarized, we define
the (total)  height $h(M)$ of $M$ by
$$h(M):=\sum_{w, d\in \Z, d\geq 0}\; h_{w,d}(M).$$
\end{para}

Like in the pure case,  the following finiteness is a basic question. We fix a type $\Psi$ of mixed motive by giving  the ranks of $\gr^W_w$ for all $w$.
\begin{conj}
Fix $c>0$ and  integers $n_w\geq 1$  ($w\in \Z$). Then there are only finitely many isomorphism classes of mixed $\Z$-motives $M$ over $K$ of type $\Psi$ with polarized graded quotients for the weight filtration, with semi-stable reductions,  such that $h(M)<c$ and such that the degree of the polarization is $(n_w)_w$.

\end{conj}
\begin{para} If we assume this conjecture, the finite generations of the motivic cohomology groups $Ext^i(\Z, M)$ ($i\geq 0$) associated to $\Z$-motives $M$ over $K$ (or their subgroups, for example  $(O_K)^\times \subset K^\times = Ext^1(\Z, \Z(1))$) are reduced by the standard arguments to the weak Mordell-Weil for them.

\end{para}

 \begin{para}\label{GA1} Let $K$ be a number field and let $C$ be a proper smooth curve over $K$. Let $M$ be a mixed $\Z$- motive over the function field $K(C)$ with polarized graded quotients for the weight filtration. Fix $w\in \Z$ and $d\geq 0$. Assume that either the following (i) or (ii) is satisfied. (i) $d\leq 3$. (ii) 
 For any  $r\in \Z$ such that $0\leq r\leq -d$, the local monodromies of $\gr^W_r(M)$ are trivial.

 Concerning the relation between geometric heights and arithmetic heights, the following formula was proved in many cases in the joint work of Spencer Bloch and the author:
$$h_{w,d}(M(t))/[L:K] = h_{w,d}(M) h(t) + O(1).$$
Here $L$ is a finite extension of $K$, $t$ is an $L$-rational point of $C$ at which $M$ does not degenerate, $h_{w,d}(M(t))$ is the height of the motive $M(t)$ over $L$ defined in section \ref{sec.ari}, $h_{w,d}(M)$ is the geometric height of $M$ defined in section \ref{sec.geo}, and $h(t)$ is the height of the point $t$ defined as $h_{\cL}(t)/\text{deg}(\cL)$, where $\cL$ is an ample line bundle $\cL$ on $C$ and $h_{\cL}$ is a height function associated to $\cL$. The $O(1)$ is a function in $t$ which is  bounded below and above independently of $L$ and $t$. 

For example, let $E$ be an elliptic curve over $K(C)$, let $a\in E(K(C))$, and let $M$ be the mixed motive over $K(C)$ corresponding to $a$ such that $W_0M=M$, 
$\gr^W_0=\Z$, $\gr^W_{-1}$ is the $H_1$ of $E$, and $W_{-2}M=0$. Then $h_{0,-1}(M(t))$ coincides with the N\'eron-Tate height of $a(t)$ and $h_{0,-1}(M)$ coincides with the geometric height of $a$ defined by using the intersection theory. The above formula for $h_{0,-1}$ in this case was proved by Tate in  \cite{Ta}. A generalization to abelian varieties was obtained in Green \cite{Gre}. 

\end{para}

\begin{para}\label{GA2} Assume furthermore $d\geq 2$. There is a related local version of \ref{GA1} which Spencer Bloch and the author proved in many cases: $$h_{w,d,v}(M(t)) = h_{w,d, x}(M) h_{x,v}(t) + O(1).$$
Here $x\in C(K_v)$ is fixed, $t\in C(K_v)$, $t\neq x$,  $t$ converges to $x$ in $C(K_v)$, $h_{w,d,v}(M(t))$ is the local height of $M(t)$ at $v$, $h_{w,d, x}(M)$ is the local height of $M$ at $x$, and $h_{x,v}(t)$ is the local height function defined as follows. By taking a local coordinate function $q$ on $C$ at $x$, it is defined as
$$h_{x,v}(t) = - \log(|q(t)|_v).$$

\end{para}

\begin{para} If neither the condition (i) nor (ii) in \ref{GA1} is satisfied, it seems that formulas like  \ref{GA1} and \ref{GA2} become more complicated.

\end{para}

Kazuya KATO,

Department of mathematics,
University of Chicago,
Chicago, Illinois, 60637,
USA,

{\tt kkato@math.uchicago.edu}

\medskip

Key words. Height, motive, Hodge theory, local monodromy, weight filtration

Running title. Heights of mixed motives


\begin{thebibliography}{99}

\bibitem{Be}
{\sc Beilinson, A.}, 
{\em  Height pairing between algebraic cycles},  
Contemp. Math. {\bf 67}, (1987), 1--24.
\bibitem{Bl}
{\sc Bloch, S.},
{\em Height pairings for algebraic cycles}, 
J. Pure Appl. Algebra {\bf 34} (1984), 119--145. 

\bibitem{CKS}
{\sc Cattani, ~E., Kaplan, A., Schmid, ~W.},
{\em Degeneration of Hodge structures},
Ann. of Math.
{\bf 123} (1986), 457--535.



\bibitem{De}
{\sc Deligne  ~P.},
{\em La conjecture de Weil},
Publ. Math. I.H.E.S. {\bf 52}
(1980), 137--252.

\bibitem{Gre}{\sc Green, W.},
{\em Heights in families of abelian varieties},
 Duke Math. J. {\bf 58} (1989), 617--632. 

\bibitem{Gri}
{\sc Griffiths, P.}, 
{\em Periods of integrals on algebraic manifolds. III. Some global differential-geometric properties of the period mapping}, 
IHES. Publ. Math. {\bf 38} (1970),  125--180. 

\bibitem{KM}
{\sc Kashiwara, ~M.},
{\em  A study of variation of mixed Hodge structure}, 
Publ. R.I.M.S., Kyoto Univ. {\bf 22} (1986),  991--1024.


\bibitem{KK}
{\sc Kato, ~K.},
{\em Heights of motives}, preprint (2013), arXiv:1306.5691 [math.NT].

\bibitem{KNU}
{\sc Kato, K., Nakayama, C., Usui, S.},
{\em Classifying spaces of degenerating mixed Hodge structures},
I., Adv. Stud. Pure Math. {\bf 54}, 187--222, II., Kyoto J. Math. {\bf 51}, 
149--261, III., to appear in J. Algebraic Geometry. 


\bibitem{Pe}
{\sc Peters, C.},
{\em A criterion for flatness of Hodge bundles over curves and geometric applications},Math. Ann. {\bf 268} (1984), 1--19.

\bibitem{Pe2}
{\sc Peters, C.},
{\em Arakelov-type inequalities for Hodge bundles}, preprint (2000), arXiv:math/0007102.


\bibitem{Sc}
{\sc Schwarz, C.}, 
{\em Relative monodromy weight filtrations},
 Math. Z. {\bf 236} (2001), 11--21.
 
 
\bibitem{Ta}{\sc Tate, J.},
\newblock
{\em Variation of the canonical height of a point depending on a parameter}
 Amer. J. Math. 105 (1983), 287--294.

\bibitem{VZ}{\sc Viehweg, E., Zuo, K}
{\em  Families over curves with a strictly maximal Higgs field}, preprint (2003), arXiv:math/0303103.
\end{thebibliography}
\end{document}